\documentclass[11pt]{amsart}
\usepackage[dvips]{graphicx}
\usepackage[all]{xy}

\newtheorem{theorem}{Theorem}

\theoremstyle{definition}

\numberwithin{equation}{section}

\theoremstyle{corollary}

\numberwithin{equation}{section}
\theoremstyle{example}

\theoremstyle{proposition}

\newcommand\R{\mathbb{R}}

\newcommand\C{\mathbb{C}}

\renewcommand\P{\mathbb{P}}

\newcommand\SU{{\rm SU}}
\newcommand\SL{{\rm SL}}

\newcommand\U{{\rm U}}

\begin{document}
%
%
\title{A Hamiltonian stable Minimal Lagrangian submanifold of projective
space with non-parallel second fundamental form}
\author{Lucio Bedulli and Anna Gori}
\address{Dipartimento di Matematica - Universit\`a di Bologna\\
Piazza di Porta S. Donato 5,
40126 Bologna\\Italy}
\email{bedulli@math.unifi.it}
\address{Dipartimento di Matematica e Appl.\ per l'Architettura -Universit\`a di Firenze\\
Piazza Ghiberti 27\\50100 Firenze\\Italy}
\email{gori@math.unifi.it}
\begin{abstract}
In this note we show that Hamiltonian stable minimal Lagrangian submanifolds of projective
space need not have parallel second fundamental form.
\end{abstract}
\maketitle
\section*{Introduction}
%
%
%
%
We describe an example of a Hamiltonian stable minimal Lagrangian
submanifold of $\C\P^3$ with non parallel second fundamental
form.\\ Let $i:L\hookrightarrow M$ be a minimal Lagrangian
submanifold of a K\"ahler manifold $(M,g,\omega)$. A deformation
$\{L_t\}$ of $L$ with $\frac{dL_t}{dt}|_{t=0}=V$ is called a {\em
Hamiltonian variation} if the $1$-form on $L$, $i^*(\iota_V
\omega)$ is exact. We call a minimal Lagrangian submanifold {\em
Hamiltonian stable} if the second variation of the area functional
through Hamiltonian variations is non negative.\\
In \cite{O}  Oh characterizes the Hamiltonian stability of minimal
Lagrangian submanifolds of K\"ahler-Einstein manifolds.\\

{\bf{Theorem}} [Oh] {\em Let $(M,\omega)$ be a K\"ahler-Einstein manifold, with Einstein constant $\kappa$.
Assume that $L$ is a minimal Lagrangian submanifold. Let
$\lambda_1(L)$ be the first eigenvalue of the Laplacian
acting on $C^{\infty}(L)$. Then $L$ is stable under Hamiltonian
deformations if and only if $\lambda_1(L)\geq \kappa$.}\\

Note that Ono in \cite{Ono} proves that, if $M$ is a Hermitian
symmetric space, a compact minimal Lagrangian submanifold $L$
actually has $\lambda_1(L)\leq \kappa$; therefore it is Hamiltonian
stable if and only if $\lambda_1(L)=\kappa$.  In particular, when
$M$ is the complex projective space $\C\P^n$ endowed with the Fubini
Study metric of constant holomorphic
sectional curvature $c$, $\lambda_1(L)$ must be equal to $c\frac{(n+1)}{2}$.\\

Amarzaya and Ohnita in \cite{AO} have studied the Hamiltonian
stability  of certain Lagrangian submanifolds of $\C\P^n$; in
particular they have shown that all the minimal Lagrangian
submanifolds of the complex projective space with parallel second
fundamental form are Hamiltonian stable. This leads them to
formulate the following question (Problem 4.2 \cite{AO}): should a
compact minimal Hamiltonian stable Lagrangian submanifold of
$\C\P^n$ have parallel second fundamental form? Our result gives a
negative answer to this question, proving
\begin{theorem} The complex projective space, $\C\P^3$, endowed with the
Fubini Study metric, admits a homogeneous minimal Lagrangian Hamiltonian
stable submanifold $L$, whose second fundamental form is not
parallel.
\end{theorem}
In order to show that $L$ is Hamiltonian stable we will perform an
explicit computation of the first eigenvalue $\lambda_1(L)$ of the
Laplacian based on general theory of invariant differential
operators on homogeneous spaces, for which we refer to \cite{T},
\cite{H}.

\section{The Lagrangian submanifold $L\subset \C\P^3$ and its least positive eigenvalue $\lambda_1(L)$}

\paragraph{\bf The Lagrangian $\SU(2)$-orbit of $\C\P^3$}
Consider the irreducible linear representation of $\SU(2)$ on
$S^3\C^2$. From now on we will regard the representation space as
the space of complex homogeneous polynomials of degree 3 in two
variables $z_1, z_2$. As proven in \cite{SK} and \cite{BGLag},
$\P(S^3\C^2)\cong \C\P^3$ is almost homogeneous under the
complexified action of $\SL(2,\C)$ and the open orbit is a Stein
manifold. From Theorem 1 in \cite {BGLag} it follows that $\C\P^3$
admits a unique $\SU(2)$-homogeneous Lagrangian submanifold $L$
which is also minimal with respect to the K\"ahler structure given
by Fubini Study metric $g_{FS}$. This submanifold turns out to be
the $\SU(2)$-orbit through the point $p=[z_1^3+z_2^3]$ and has
already been studied, from the topological viewpoint, by Chiang in
\cite{Ch}.\\
The second fundamental form of the Lagrangian orbit is not parallel
(in \cite{N} the complete classification of compact Lagrangian
submanifolds of $\C\P^n$ with parallel second fundamental form is
given).
\paragraph{\bf The isotropy representation}
The isotropy at the point $p$ is given by the finite
subgroup $F$ of $\SU(2)$ generated by
\[
  a=\left[
\begin{array}{cc}
e^{i\frac{\pi}{3}} & 0 \\
  0 & e^{-i\frac{\pi}{3}}
\end{array}
\right],\: \quad b= \left[
\begin{array}{cc}
0 & i \\
i & 0
\end{array}
\right] \:.
\]
The isotropy representation splits the tangent space to this Lagrangian orbit at $p$ as the direct sum
of two $F$-invariant irreducible submoduli $V_1\oplus V_2$, where the first summand is 1-dimensional and
$F$ acts as ${-Id}$ on it.\\
Our argument in this computation is analogous to the one used by
Muto and Urakawa in \cite{UM}.
They assume the existence of a 1-dimensional invariant subspace
fixed under the action of the isotropy subgroup. In our context
there is no $F$-fixed subspace, nevertheless we will see that for
our prupose it is sufficient to have a vector $v$ such that
$g^2\cdot v=v$ for all $g \in F$. This is true for $v \in V_1$.
\paragraph{\bf The induced metric on $L$}
Fix an orthonormal basis with respect to the opposite of the
Killing form $B$ on $\mathfrak {su}(2)$. Recall that if we see
$\mathfrak {su}(2)$ as a matrix group, $B$ has the following form
$$B(\xi,\eta):=4\mbox{Tr}(\xi\eta)$$
for all $\xi,\eta\in \mathfrak{su}(2)$.\\
Consider
\[ H=
\left[
\begin{array}{cc}
i & 0 \\
0 & -i
\end{array}
\right],\:\: X= \left[
\begin{array}{cc}
0 & 1 \\
-1 & 0
\end{array}
\right],  \:\:\: Y= \left[
\begin{array}{cc}
0 & i \\
i & 0
\end{array}
\right]
\]
in $\mathfrak{su}(2)$. Then $X_1:=H/2\sqrt{2},X_2:=X/2\sqrt{2}$
and $X_3:=Y/2\sqrt{2}$ form an orthonormal basis with respect to
$-B$. Denote by $\{w_j\}_{j=0\ldots 3}$ the homogeneous complex
coordinates with respect to the basis
$\{z_1^{3-j}z_2^j\}_{j=0\ldots 3}$ on $\P(S^3\C^2)$ so that, in
this frame, $p=[1:0:0:1]$. Consider around $p$ the real
coordinates $x_j:=\mbox{Re}\frac{w_j}{w_0}$ and $y_j=\mbox{Im}\frac{w_j}{w_0}$ for $j=1,2,3$.\\
The fundamental fields associated to $X_1,X_2,X_3$ at the point
$p$ then can be written as
$$
\widehat{X}_{1p}=\frac{-3}{\sqrt{2}}\frac{\partial}{\partial y_3},
$$
$$
\widehat{X}_{2p}=\frac{-3}{2\sqrt{2}}\frac{\partial}{\partial x_1}+
\frac{3}{2\sqrt{2}}\frac{\partial}{\partial x_2},
$$
$$
\widehat{X}_{3p}=\frac{3}{2\sqrt{2}}\frac{\partial}{\partial y_1}+
\frac{3}{2\sqrt{2}}\frac{\partial}{\partial y_2},
$$
where we drop
the subscript $p$ from the basis on the tangent space. These three
vectors span the tangent space to the orbit $L$ at $p$. Note that $V_1$ is generated by
$\widehat{X}_{1p}$ while $V_2$ is generated by $\widehat{X}_{2p}$ and $\widehat{X}_{3p}$. Now we
explicitly compute the metric $g$ induced by $g_{FS}$ on $T_p L$ with respect
to the basis $\{\widehat{X}_{1p},
\widehat{X}_{2p},\widehat{X}_{3p}\}$. To this aim we follow the
general construction of the Fubini Study metric of constant
holomorphic sectional curvature $c=4$ on $\C\P^n$ (see e.g.
\cite{KN} p. 273), starting from a Hermitian metric $h$ on
$\C^{n+1}$. This makes the fibration of the unit sphere of
$(\C^{n+1},h)$ onto $\C\P^n$ a Riemannian submersion.
An $\SU(2)$-invariant Hermitian metric on $S^3\C^2$ is of the form
$$h=u\;dw_1\otimes d\overline{w}_1+\frac{1}{3}u\;dw_2\otimes d\overline{w}_2+\frac{1}{3}u\;
dw_3\otimes d\overline{w}_3+ u\;dw_4\otimes d\overline{w}_4$$ with
$u>0$.
The horizontal lifting at
$\widetilde{p}=(\frac{1}{\sqrt{2u}},0,0,\frac{1}{\sqrt{2u}})\in
S^7$ of a generic tangent vector
$$
v=\sum_{j=1}^3 a_j\frac{\partial}{\partial x_j}+\sum_{j=1}^3 b_j\frac{\partial}{\partial
 y_j}
$$
is then given by
$$
\frac{1}{2\sqrt{2u}}(-a_3\frac{\partial}{\partial \widetilde{x}_0}-b_3\frac{\partial}
{\partial \widetilde{y}_0}+2a_1\frac{\partial}{\partial \widetilde{x}_1}+2b_1\frac{\partial}
{\partial \widetilde{y}_1}+2a_2\frac{\partial}{\partial \widetilde{x}_2}+2b_2\frac{\partial}
{\partial \widetilde{y}_2}+a_3\frac{\partial}{\partial \widetilde{x}_3}+
b_3\frac{\partial}{\partial \widetilde{y}_3}),
$$
where $\{\widetilde{x}_j,\widetilde{y}_j\}_{j=0\ldots 3}$ are the real coordinates of $\C^4$
corresponding to $\{w_j\}$. Therefore
$$
\|\widehat{X}_{1p}\|^2_g=\frac{9}{8},\quad \quad
\|\widehat{X}_{2p}\|^2_g=\frac{3}{8},\quad \quad
\|\widehat{X}_{3p}\|^2_g=\frac{3}{8}.
$$
Hence
$\frac{2\sqrt{2}}{3}\widehat{X}_{1p}$,
$\frac{2\sqrt{2}}{\sqrt{3}}\widehat{X}_{2p}$ and
$\frac{2\sqrt{2}}{\sqrt{3}}\widehat{X}_{3p}$ form a
$g$-orthonormal basis for $T_p L$.\\

\paragraph{\bf The computation of $\lambda_1(L)$}
First recall well known results that can be found e.g. in \cite{T}
\cite{H}. Let $M=G/K$ be a $n$-dimensional homogeneous space,
where $G$ is a compact connected Lie group and $K$ a closed
subgroup of $G$. The choice of an $Ad(G)$-invariant scalar product
on the Lie algebra $\mathfrak{g}$ gives an orthogonal splitting
$\mathfrak{g}=\mathfrak{k}\oplus\mathfrak{m}$. Let
$S(\mathfrak{m})_K$ be the algebra of $Ad(K)$-invariant elements
of the symmetric algebra of $\mathfrak{m}$. Denote by
$S(\mathfrak{m})_K^\C$ its complexification. Fix a basis $\{Y_i\}$
of $\mathfrak{m}$ and regard the elements $P$ of
$S(\mathfrak{m})_K^\C$ as polynomials in $Y_1,\ldots,Y_n$.\\
Denote by $\mathcal{D}(M)$ the space of all $G$-invariant
differential operators on $C^\infty(M;\C)$.
\begin{theorem}
The map $\widehat\lambda \colon S(\mathfrak{m})_K^\C \to {\mathcal{D}}(M)$
defined by
\[
[\widehat\lambda(P(Y_1,\ldots,Y_n))f](x\cdot
K)=[P(\frac{\partial}{\partial y_1},\cdots,\frac{\partial}{\partial
  y_n})f(x\: {\rm exp}(\sum y_iY_i)\cdot K)](0),
\]
is a linear isomorphism. Moreover if $\{Y_1,Y_2,\cdots,Y_n\}$ is an
orthonormal basis with respect to an inner product on
$\mathfrak{m}$ and $\Delta_g$ is the corresponding Laplacian of
the $G$-invariant metric $g$ induced on $M$, then
\[
\widehat\lambda(\sum_iY_i^2)=-\Delta_g.
\]
\end{theorem}

Fix a unitary representation $\rho \colon G \to
U(V_{\rho},\langle\,,\rangle)$ of dimension $d_\rho$. Denote by $V_\rho^K$
the set of vectors fixed by $K$ and by $m_\rho$ its dimension. Let
$\{v_i\}_{i=1}^{d_\rho}$ be an orthonormal basis of $V_{\rho}$ whose
first $m_\rho$ vectors span $V_\rho^K$. Define the complex valued
functions $\rho_{ij}(xK)=\langle \rho(x)v_j,v_i\rangle$ for
$i=1,\ldots,d_\rho$ and
$j=1,\ldots,m_\rho$.
A representation $\rho$ such that $m_{\rho}>0$ is said to be a
{\em spherical representation} for the pair $(G,K)$.\\
Peter-Weyl Theorem (see e.g. \cite{T}, \cite{H}) says that
$\{\sqrt{d_{\rho}}\overline{\rho}_{ij}\}$ is a complete orthonormal system of $L^2(M,\C)$
with respect to the standard $L^2$-norm corresponding to the
$G$-invariant Riemannian metric $g$ on $M$ induced by an
$Ad(K)$-invariant inner product on
$\mathfrak{m}$.\\
It is possible to prove that $\overline{\rho}_{ij}$ are eigenfunctions for the Laplacian
$\Delta_g$ and  ``Freudenthal Formula'' gives us the eigenvalues:
\[
\Delta_g \overline{\rho}_{ij}=(\mu_\rho+2\delta,\mu_\rho)\overline{\rho}_{ij},
\]
where $\delta=\frac{1}{2}\sum_{\alpha\in \Delta^+}\alpha$,
$\Delta^+$ is the set of positive roots,  $\mu_\rho$ is the highest
weight of $\rho$ and $(\:\:,\:\:)$ is the inner product on
$\mathfrak{g}^*$ induced by the Killing
form.\\
\\
We return now to the homogeneous Lagrangian submanifold
$\SU(2)/F$.\\ The basis in $\mathfrak{m}$ corresponding to
$\{\frac{2\sqrt{2}}{3}\widehat{X}_{1p},
\frac{2\sqrt{2}}{\sqrt{3}}\widehat{X}_{2p},
\frac{2\sqrt{2}}{\sqrt{3}}\widehat{X}_{3p}\}$ will be denoted by
$\{Y_1,Y_2,Y_3\}$. The isotropy $F$ at $p$ acts on the
$1$-dimensional space spanned by $\widehat{X}_{1p}$ as $-Id$ so
that
\begin{equation}
\widehat{\lambda}(X_1^2) \overline{\rho}_{ij} (x\cdot
F)=\overline{\langle\rho(x)d\rho(X_1)^2v_j,v_i\rangle},
\end{equation}
for every $x\in \SU(2)$, $\rho$ representation of $\SU(2)$
($d\rho$ is the representation $\rho$ at the Lie algebra level).\\
 Now
$$
\Delta_g=-\widehat{\lambda}(Y_1^2+Y_2^2+Y_3^2)=-\widehat{\lambda}
(\frac{8}{9}X_1^2+\frac{8}{3}X_2^2+\frac{8}{3}X_3^2)
=-\frac{8}{3}\widehat{\lambda}({X}_1^2+{X}_2^2+{X}_3^2)-\frac{16}{9}\widehat{\lambda}(-{X}_1^2),
$$
where the last equality makes sense for ${X}_1^2\in
S(\mathfrak{m})_F.$ Hence
$$\Delta_g\overline{\rho}_{ij}=\frac{8}{3}\Delta_{-B}-
\widehat{\lambda}(-\frac{16}{9}{X}_1^2)\overline{\rho}_{ij}=\frac{8}
{3}(\mu_{\rho}+\alpha^+,\mu_{\rho})\overline{\rho}_{ij}-\widehat{\lambda}
(-\frac{2}{9}H^2)\overline{\rho}_{ij}.$$ So we have to study the
operator ${D}=d\rho(H)^2$.\\
The following remark is analogous to Lemma 3.1
in \cite{UM}. For every spherical unitary representation
$\rho:\SU(2)\to \U(V_{\rho},\langle\,,\rangle)$ for the pair
$(\SU(2),F)$ the operator $ {D}:V_{\rho}\to V_{\rho}$ is
self-adjoint. This follows from the $\SU(2)$-invariance of
$\langle\,,\rangle$; indeed this implies
\[
\langle d\rho(H)u,v\rangle+\langle u,d\rho(H) v\rangle=0
\]
for every $u,v\in V_{\rho}$.\\
Moreover the subspace $V_{\rho}^F$ is ${D}$-invariant, in fact for every
$v \in V_\rho^F$, since $Ad(a)H=H$ and $Ad(b)H=-H$, we have
\begin{eqnarray*}
\rho(a)({D}v) & = &\rho(a)({d\rho(H)^2}v)=\rho(a)({d\rho(Ad(a)H)^2}v)\\
 & = & \rho(a)\rho(a^{-1})d\rho(H)\rho(a)\rho(a^{-1})d\rho(H)\rho(a)v={D}v
\end{eqnarray*}
and
\begin{eqnarray*}
\rho(b)({D}v) & = & \rho(b)({d\rho(H)^2}v) =\rho(b)(d\rho(Ad(b)H)^2v)\\
 & = & \rho(b) (-d\rho(Ad(b)H))(-d\rho(Ad(b)H))v={D}v.
\end{eqnarray*}
Then we can find an orthonormal basis
$\{u_j\}_{j=1}^{d_{\rho}}$ for $V_{\rho}$, of eigenvectors for
$D$ so that the first $m_{\rho}$ vectors are a basis for $V_{\rho}^F$.
Denote by $\mu_j$ the corresponding eigenvalues.
In terms of this basis
\[
\widehat{\lambda}(H^2)\overline{\rho}_{ij}(x\cdot F
)=\overline{(\rho(x)\cdot d\rho(H)^2 u_j,u_i)}=\mu_j
\overline{\rho}_{ij}(x\cdot F)
\]
and
$$\Delta_g\overline{\rho}_{ij}=\frac{8}
{3}(\mu_{\rho}+\alpha^+,\mu_{\rho})\overline{\rho}_{ij}+\frac{2}{9}
\mu_j\overline{\rho}_{ij}.$$
Thus we have to find
\begin{equation}\label{minaut}\lambda_1(L)=\min_\rho
\{\lambda_1^{\rho}\},\quad\mbox{where}\quad
\lambda_1^{\rho}=\min_{1\leq j\leq m_{\rho}}\{\frac{8}
{3}(\mu_{\rho}+\alpha^+,\mu_{\rho})+\frac{2}{9}\mu_j\}
\end{equation}
as $\rho$ varies through the irreducible spherical representations
of the pair $(\SU(2),F)$. Recall
that all the irreducible representations of $\SU(2)$ are the
symmetric powers of the standard representation of $\SU(2)$. If
$\rho=\rho_k=S^k(\C^2)$, then $\mu_{\rho}:\mathfrak{t}\to i\R$ is
defined by $\mu_{\rho}(H)=ik$. Hence, since the positive
root $\alpha^+$ sends $H$ to $2i$, we have $(\mu_{\rho}+\alpha^+,\mu_{\rho})=\frac{1}{4}k+\frac{1}{8}k^2$.\\
Now we determine the subspaces of $V_{\rho_k}$ fixed by the isotropy $F$.
We will see $V_{\rho_k}$ as the space of complex homogeneous polynomials of degree $k$ in
two variables $z_1,z_2$, so that we can write a generic vector $v$ as
$$
\sum_{j=0}^k c_j z_1^j z_2^{k-j}.
$$
Now $bv=i^k\sum_{j=0}^k c_j z_2^j z_1^{k-j}$. The condition $bv=v$
forces $k$ to be even and imposes the following constraints:
\\
$\bullet$ If $k\equiv_4 0$ then, for every $j$, $c_j=c_{k-j}$; \\
$\bullet$ If $k\equiv_4 2$ then, for every $j$, $c_j=-c_{k-j}$.\\
Furthermore the condition $av=\sum_{j=0}^k c_j \alpha^{2j-k}z_1^j
z_2^{k-j}=v$ leads to distinguish the following cases, where
$[\![\,,\,]\!]$ stands for ``the span of'' \\
$\bullet$ if $k\equiv_{12} 0$ then $V^F_{\rho_k}=[\![z_1^k+z_2^k,
z_1^{k-3}z_2^3+z_2^{k-3}z_1^3,\ldots,
z_1^{\frac{k}{2}}z_2^{\frac{k}{2}}]\!];$
\\
$\bullet$ if $k\equiv_{12} 2$ then
$V^F_{\rho_k}=[\![z_1^{k-1}z_2-z_2^{k-1}z_1,
z_1^{k-4}z_2^4-z_2^{k-4}z_1^4, \ldots,
z_1^{\frac{k}{2}+1}z_2^{\frac{k}{2}-1}-z_1^{\frac{k}{2}-1}z_2^{\frac{k}{2}+1}]\!]$;
\\
$\bullet$ if $k\equiv_{12} 4$ then $V^F_{\rho_k}=[\![
z_1^{k-2}z_2^2+z_2^{k-2}z_1^2,
z_1^{k-5}z_2^5+z_2^{k-5}z_1^5,\ldots,
z_1^{\frac{k}{2}+2}z_2^{\frac{k}{2}-2}+z_1^{\frac{k}{2}-2}
z_2^{\frac{k}{2}+2}]\!];$
\\
$\bullet$ if $k\equiv_{12} 6$ then $V^F_{\rho_k}=[\! [z_1^k-z_2^k,
z_1^{k-3}z_2^3-z_2^{k-3}z_1^3,\ldots,
z_1^{\frac{k}{2}+3}z_2^{\frac{k}{2}-3}-z_1^{\frac{k}{2}-3}
z_2^{\frac{k}{2}+3}]\!];$
\\
$\bullet$ if $k\equiv_{12} 8$ then $V^F_{\rho_k}=[\![
z_1^{k-1}z_2+z_2^{k-1}z_1, z_1^{k-4}z_2^4+z_2^{k-4}z_1^4, \ldots,
z_1^{\frac{k}{2}+1}z_2^{\frac{k}{2}-1}+z_1^{\frac{k}{2}-1}
z_2^{\frac{k}{2}+1}]\!];$
\\
$\bullet$ if $k\equiv_{12} 10$ then $V^F_{\rho_k}=[\![
z_1^{k-2}z_2^2-z_2^{k-2}z_1^2, z_1^{k-5}z_2^5-z_2^{k-5}z_1^5,
\ldots, z_1^{\frac{k}{2}+2}z_2^{\frac{k}{2}-2}-z_1^{\frac{k}{2}-2}
z_2^{\frac{k}{2}+2}]\!].$
\\
Actually the generators of $V^F_{\rho_k}$ we wrote down are
eigenvectors for  $d\rho(H)^2$. Indeed this follows from the fact
that
$$
d\rho(H)z_1^lz_2^{k-l}=\frac{d}{dt}_{|_t=0}\exp tH\cdot z_1^lz_2^{k-l}=i(2l-k)z_1^lz_2^{k-l}
$$
and
$$
d\rho(H)^2 z_1^lz_2^{k-l}=-(2l-k)^2 z_1^lz_2^{k-l},
$$
where \\
$\bullet$ for $k\equiv_{12} 0$, $l=k,k-3,\dots,\frac{k}{2}$;
\\
$\bullet$ for $k\equiv_{12} 2$, $l=k-1,k-4,\dots,\frac{k}{2}+1$;
\\
$\bullet$ for $k\equiv_{12} 4$, $l=k-2,k-5,\dots,\frac{k}{2}+2$;
\\
$\bullet$ for $k\equiv_{12} 6$, $l=k,k-3,\dots,\frac{k}{2}+3$;
\\
$\bullet$ for $k\equiv_{12} 8$, $l=k-1,k-4,\dots,\frac{k}{2}+1$;
\\
$\bullet$ for $k\equiv_{12} 10$, $l=k-2,k-5,\dots,\frac{k}{2}+2$.
\\ \\
For every $k$,  $\lambda_1^{\rho_k}$ in (\ref{minaut}) is attained when
$l=l_{\mbox{m}}$ is maximal.
Therefore, when $k\equiv_6 0$ then
$l_{\mbox{m}}=k$ and $\lambda_1^{\rho_k}=\frac{2}{3}k+\frac{1}{9}k^2$,
when $k \equiv_6 2$ ($k \geq 6$) then $l_{\mbox{m}}=k-1$ and
$\lambda_1^{\rho_k}=\frac{1}{9}(k^2+14k-8)$, while when $k\equiv_6 4$
then $l_{\mbox{m}}=k-2$ and $\lambda_1^{\rho_k}=\frac{1}{9}(k^2+22k-32).$
Finally the non-zero eigenvalues are always greater or equal than
8, and the equality holds in the first and third case for $k=6$ and $k=4$ respectively.
 Our claim then follows since in these two cases
$\lambda_1(L)$ is equal to the Einstein constant $\kappa=2c$.

\end{document}